\documentclass[amscd,amssymb,verbatim,12pt]{amsart}
\newif\ifAMS
\IfFileExists{amssymb.sty}
  {\AMStrue\usepackage{amssymb}}
  {\usepackage{latexsym}}
  \usepackage{enumerate}
\theoremstyle{plain}

\newtheorem*{Gromov}{Gromov's Theorem}
\newtheorem*{sphere}{Proposition \ref{P:sphere} }
\newtheorem*{cheeger}{Theorem \ref{T:cheeger}}
\newtheorem*{profile}{Theorem \ref{T:profile1}}
\newtheorem*{profcor}{Corollary \ref{C:profile2}}
\newtheorem*{filling}{Theorem \ref{T:dim3}}
\newtheorem*{Corgro}{Corollary}

\newtheorem{Thm}{Theorem}[section]
\newtheorem{Cor}[Thm]{Corollary}
\newtheorem{Lem}[Thm]{Lemma}
\newtheorem{Prop}[Thm]{Proposition}

\theoremstyle{definition}
\newtheorem{Def}{Definition}

\theoremstyle{remark}

\newtheorem{Rem}{Remark}

\newtheorem{Ex}[Thm]{Example}

\newcommand{\interior}{^{ \kern-5pt ^\circ}}
\newcommand {\bd}{\partial}

\begin{document}
\title
{Cheeger constants of surfaces and isoperimetric inequalities}

\author
{Panos Papasoglu }


\email {panos@math.uoa.gr}

\address
{Mathematics Department, University of Athens, Athens 157 84,
Greece }
\begin{abstract}
We show that the Cheeger constant of compact surfaces is bounded
by a function of the area. We apply this to isoperimetric profiles
of bounded genus non-compact surfaces, to show that if their
isoperimetric profile grows faster than $\sqrt t$, then it grows
at least as fast as a linear function. This generalizes a result
of Gromov for simply connected surfaces.

We study the isoperimetric problem in dimension 3. We show that if
the filling volume function in dimension 2 is Euclidean, while in
dimension 3 is sub-Euclidean and there is a $g$ such that
minimizers in dimension 3 have genus at most $g$, then the filling
function in dimension 3 is `almost' linear.

\end{abstract}
\maketitle
\section{Introduction}

If $M$ is a riemannian manifold of dimension $n$ one defines the
Cheeger constant $h$ of $M$ by
$$h(M)=\underset {A} { \inf } \, \{\frac {vol_{n-1}(\partial A)}
{vol_n(A)}: vol_n(A)\leq \frac {1} {2} vol_n(M) \}$$ where $A$
ranges over all open subsets of $M$ with smooth boundary. If $M$
is a simplicial manifold one can define the Cheeger constant of
$M$ similarly.


As usual we call $vol_2$ area and $vol_1$ length. If $M$ is a
simplicial 2-manifold or a 2-manifold with a riemannian metric we
denote by $A(M)$ the area of $M$. Similarly if $p$ is a
(simplicial or riemannian) path we denote by $l(p)$ the length of
$p$.


We will show that one can give a bound of the Cheeger constant of
a surface that depends only on its area. So for example there is a
constant $c$ such that any riemannian manifold homeomorphic to the
2-sphere $S$, which has area 1, has $h(S)\leq c$. We state our
results both in the simplicial and in the riemannian setting. Our
results in the simplicial case are applied  in the last section to
higher isoperimetric inequalities. We provide explicit bounds but
the constants in the theorems below are far from optimal.

\begin{sphere}  Let $S$ be a riemannian manifold
or a simplicial complex homeomorphic to the 2-sphere. Then the
Cheeger constant, $h(S)$, of $S$ satisfies the inequality:
$$ h(S)\leq \frac {16} {\sqrt{A(S)}} $$
where $A(S)$ is the area of $S$.
\end{sphere}

In general we obtain an upper bound that depends on the genus:

\begin{cheeger} Let $S$ be a closed orientable
surface of genus $g\geq 1$ equipped either with a riemannian
metric or with a simplicial complex structure. Let $A(S)$ be its
(simplicial or riemannian) area. Then the Cheeger constant,
$h(S)$, of $S$ satisfies the inequality:
$$ h(S)\leq \frac {4\cdot 10^3\cdot g^2} {\sqrt{A(S)}} $$
\end{cheeger}

One sees easily that the same bound applies to surfaces with
boundary (just collapse the boundary curves to points to obtain a
closed surface). One can get bounds for non-orientable surfaces
too by passing to the orientable double cover.

If $(M^n,g)$ is a riemannian manifold of infinite volume the
isoperimetric profile function of $M^n$ is a function $I_M:\Bbb R
^+\to \Bbb R ^+$ defined by:
$$ I_M(t)=\underset {\Omega} { \inf } \{vol_{n-1}(\bd \Omega):
\Omega\subset M^n,\, vol_n(\Omega)=t \}$$ where $\Omega $ ranges
over all regions of $M^n$ with smooth boundary. One can define
similarly an isoperimetric profile function $I_M:\Bbb N\to \Bbb N$
for simplicial manifolds $M^n$.

Other functions related to the isoperimetric problem are the
filling area, $FA_0 $, and homological filling area, $FA^h $,
functions of $M$ that we define below. For more information on
filling invariants and applications we refer the reader to the
seminal paper of Gromov \cite{Gr1}.

If $p$ is a smooth contractible closed curve in $M$ we define its
filling area, $FillA_0(p)$, as follows: We consider all riemannian
discs $D$ such that there is a 1-lipschitz map $f:D\to X$ with
$f|_{\bd D}=p$. We define $FillA_0(p)$ to be the infimum of the
areas of this collection of disks. We define now the filling area
function of $M$ by:

$$ FA_0(t)=\underset {p} { \sup } \{FillA_0(p): l(p)\leq t \}$$

where $p$ ranges over all smooth contractible closed curves of $M$
and $D$ over riemannian disks filling $p$.

 More generally we can
consider 1-cycles $c$ (i.e. unions of closed curves) that can be
filled by 2-cycles to define the homological filling area function
(see sec. 2 for details).

%

Gromov (\cite {Gr2}, ch. 6, see also \cite{CDP}, ch.6) showed the
following:
\begin{Gromov}\label{T:gromov}  Let $(M^n,g)$ be a simply connected
riemannian manifold. Assume that there is some $t_0$ such that for
all $t>t_0$, $ FA_0(t)\leq \frac {1}{16\pi }t^2$. Then there is a
constant $K$ such that for all $t>t_0$, $ FA_0(t)\leq Kt$.
\end{Gromov}

Gersten \cite {Ge} observed that this theorem holds also for
homological filling area $FA^h$ (see also \cite {Gr2}, 6.6E,
6.6F), while Olshanskii \cite{O} gave an elementary proof of
Gromov's theorem (see as well \cite{Bo}, \cite {P}, \cite {D}, for
other proofs).

If the dimension of $M$ is 2 then there is an obvious link between
filling area and isoperimetric profile, so from Gromov's theorem
we readily obtain the following:
\begin{Corgro} Let $(S,g)$ be a riemannian manifold homeomorphic to the plane.
Assume that there is some $t_0$ such that for all $t>t_0$, $
I_S(t)\geq 4\sqrt \pi \sqrt t$. Then there is a constant $\delta
>0 $ such that for all $t>t_0$, $ I_S(t)\geq \delta t$.
\end{Corgro}

We note that the isoperimetric problem for surfaces has been
studied extensively (see \cite{BC}, \cite {GP}, \cite{MHH},
\cite{Ri1}, \cite{Ri2}, \cite{T}).

 We see that the `gap' in the filling functions implies a `gap' for the
isoperimetric profiles of riemannian planes. It is reasonable to
ask whether there are gaps in the isoperimetric profile of other
surfaces. Although this does not hold in general we show that this
is true for planes with holes or more generally surfaces of finite
genus.

\begin{profile}  Let $S$ be a
plane with holes equipped either with a riemannian metric or with
a simplicial complex structure. Assume that there is some $K>0$
such that for all $t\in [K,100K]$, $ I_S(t)\geq 10^2 \sqrt t$.
Then there is a constant $\delta >0 $ such that for all $t>K$, $
I_S(t)\geq \delta t$.
\end{profile}

One obtains as a corollary that the same holds for finite genus
surfaces:

\begin{profcor}\label{T:profile2}  Let $S$ be a non-compact
surface of finite genus equipped either with a riemannian metric
or with a simplicial complex structure. Assume that there is some
$K>0$ such that for all $t\in [K,100K]$, $ I_S(t)\geq 10^2 \sqrt
t$. Then there is a constant $\delta >0 $ such that for all $t>K$,
$ I_S(t)\geq \delta t$.
\end{profcor}
It is an interesting question whether Gromov's theorem on filling
area has an analogue for higher dimensional filling functions. Our
results on Cheeger constants of surfaces can be used to obtain
some partial results in this direction. We will state our results
in the convenient setting of simplicial complexes. We remark that
if $M$ is a compact riemannian manifold then the filling functions
of its universal covering, $\tilde M $, are determined (up to some
lipschitz constant) by $\pi _1 (M)$ (see \cite {Ep}, theorems
10.3.3, 10.3.1 and \cite{BT}), so one can forget the riemannian
metric and work with a triangulation and simplicial chains to
calculate the filling functions of $\tilde M$.

To fix ideas when we refer to chains and cycles we mean always
chains and cycles for simplicial homology with $\Bbb Z$
coefficients (in fact our results apply to $\Bbb Z _2$
coefficients as well). We denote by $H_n(X)$ the $n$th-homology
group of the space $X$ with $\Bbb Z$-coefficients.

Let $X$ be a simplicial complex such that $H_k(X)=0$. If
$$S=\sum n_i\sigma _i$$ is a (simplicial) $k$-chain we define
the $k$-th volume of $S$ by $vol_k(S)=\sum |n_i|$. If $S$ is a
$k$-cycle we define the filling volume of $S$ by
$$Fillvol_{k+1}(S)=\inf \{ vol_{k+1}T: \bd T=S \}$$
For $k\geq 1$ we define the $(k+1)$th-filling volume function,
$FV_{k+1}$, of $X$ by:
$$FV_{k+1}(n)=\sup \{Fillvol_{k+1}(S): S \text { is a $k$-cycle such that } vol_k(S)\leq n \}$$
If $S$ is a $k$-cycle such that $ vol_k(S)\leq n$ and
$Fillvol_{k+1}(S)=FV_{k+1}(n)$ we say that $S$ is a
\textit{minimizer} for $FV_{k+1}(n)$.

If $S$ is a 2-cycle one can define the \textit{genus} of $S$.
Indeed $S$ is represented by a map $f:\Sigma \to X$ where $\Sigma
$ is a closed surface and $f$ is simplicial and 1-1 on open
2-simplices (see \cite {Ha}, sec.2.1, p.109). We define the genus
of $S$ to be the genus of the surface $\Sigma $.

We remark that $FA^h=FV_2$. As we noted earlier Gromov's theorem
applies to $FA^h$ as well so we have that if $$\lim_{n\to \infty }
\frac {FV_2(n)}{n^2}=0$$ then there is some $K>0$ such that
$FV_2(n)\leq Kn$ for all $n\in \Bbb N$. In general we say that
$FV_k$ is euclidean if there is some $K>0$ such that
$$ \frac{1}{K}n^{\frac{k}{k-1}}\leq FV_k(n)\leq Kn^{\frac{k}{k-1}}, \ \ \forall n\in \Bbb N$$
and we say that $FV_k$ is sub-euclidean if
$$\lim_{n\to \infty } \frac {FV_k(n)}{n^{\frac{k}{k-1}}}=0 $$
So by Gromov's theorem if $FV_2$ is sub-euclidean then it is
linear. We note that a naive guess that if $FV_3 $ is
sub-euclidean then it is linear is contradicted by Pansu's theorem
(\cite{Pa}), $FV_3$ in Heisenberg's group grows like $n^{\frac
{4}{3}}$.
 On the
other hand Gromov conjectures (\cite{Gr2}, sec. $6B_2$) that if
$X$ is a $Cat(0)$ space with a co-compact group action then
sub-euclidean filling implies linear filling in any dimension.
More generally it is believed that this is true for spaces
satisfying a cone-type inequality (see \cite{We2}).

 Another possible direction is to examine all filling functions
simultaneously. Some specific conjecture is: If $FV_2$ is bounded
by a quadratic function and $FV_3$ is sub-euclidean then $FV_3$ is
bounded by a linear function. Of course one can state this
conjecture in any dimension: If $FV_i$ is euclidean for
$i=2,...,k-1$ and $FV_k$ is sub-euclidean then $FV_k$ is linear.
The following theorem is giving some evidence
in favor of this conjecture.

\begin{filling} Let $X$ be a simplicial
complex such that $H_1(X)=H_2(X)=0$. Assume that the following
hold:
\begin{itemize} \item There is some $K>0$ such that
 $FV_2(n)\leq Kn^2$ for all $n\in \Bbb N$.

 \item $$\lim _{n\to \infty }\frac {FV_3(n)}{n^{\frac {3}{2}}}=0$$

 \item There is some $g\in \Bbb N$ such that if $S$ is a minimizer
 $2-cycle $ in $X^2$ then $S$ is represented by a surface of genus at
 most $g$.

\end{itemize}
Then for every $\epsilon >0$ we have
$$\lim _{n\to \infty }\frac {FV_3(n)}{n^{1+\epsilon }}=0$$
\end{filling}

\subsection {Outline of the proofs}

The proofs of proposition \ref{P:sphere} and theorem
\ref{T:cheeger} are based on Besicovitch lemma (and more generally
co-area inequalities).
The idea can be grasped easier in the case of the sphere. We
consider a minimal length simple closed curve $p$ that subdivides
the sphere $S$ into two pieces $S_1,S_2$ such that $A(S_1),A(S_2)$
have both area bigger than $A(S)/4$. Let's say $A(S_1)\geq
A(S_2)$. Now $S_1$ is a disk and subpaths of $p$ of length $<
l(p)/2$ are geodesic in $S_1$ because $p$ is minimal. Applying
Besicovitch lemma one sees that $l(p)$ is smaller than $4\sqrt
{A(S)}$. This implies proposition \ref{P:sphere}. We note that the
constant 4 is not optimal. This is also the case for all other
estimates we obtain in this paper. We just try to show the
existence of certain constants and we are always generous in our
estimates as our methods are not suited for approaching the best
values.

To prove theorem \ref{T:profile1} we show in fact that $I_S(t)\geq
\frac {1}{\sqrt K}t$ for all $t>K$. To show this we argue by
contradiction. We take a 1-cycle $c$ of minimal filling area such
that $FillA(c)>K$ and $l(c)<\frac {1}{\sqrt K}FillA(c)$. If $c=\bd
R$ then we collapse the `holes' of $R$ to points to get a
 sphere $S$. Applying proposition \ref{P:sphere} to $S$ we
find a 1-cycle $\gamma $ in $S$ such that $K\leq FillA (\gamma
)\leq 4K$ which satisfies the inequality $$FillA (\gamma ) > \frac
{l(\gamma )^2}{100}$$ This is somewhat tricky and the proof uses
also Besicovitch lemma and exploits the convexity of $x^2$. Note
though that intuitively the existence of such a curve is obvious
since $S$ is `positively curved' at some points. Finally we lift
$\gamma $ back to $R$ and we get a cycle with smaller filling area
than $c$ that has the same properties as $c$; this contradicts the
minimality of $c$. It is an easy corollary that the theorem holds
for finite genus surfaces in general.

The proof of theorem \ref{T:dim3} is by contradiction. We assume
that for some $\epsilon >0$ we have $$\lim _{n\to \infty }\frac
{FV_3(n)}{n^{1+\epsilon }}=\infty $$ We take $M$ `big enough' and
we consider a 2-cycle $S$ which is a minimal area counterexample
to $FV_3(n)\leq Mn^{1+\epsilon }$. We show that there is a $\delta
=\delta (\epsilon,g, K)>0$ such that $diam (S)>\delta \sqrt
{A(S)}$. This is where we use our results on Cheeger constants. To
simplify let's say that $S$ is a sphere. By proposition
\ref{P:sphere} there is a simple closed curve $p$ of length
smaller than $4\sqrt {A(S)}$ on $S$ which divides it into two
pieces with comparable area. Now if $diam (S)$ is small the
filling area of $p$ is much smaller than $A(S)$. So one can
subdivide $S$ into two 2-cycles of area roughly between $A(S)/4$
and $3A(S)/4$. It follows by the convexity of $n^{1+\epsilon }$
that one of these two 2-cycles is also a counterexample to
$FV_3(n)\leq Mn^{1+\epsilon }$, a contradiction. Given now that
the diameter of $S$ is big we take a minimal volume 3-cycle $R$
filling $S$ and we fix a point $x$ on $S$. We consider `balls'
around $x$ in $R$ and using again the convexity of
$n^{1+\epsilon}$ and an elementary inequality (lemma
\ref{L:outer}) we see that the 2-cycle, say $S_1$, given by the
boundary of some of these `balls' has filling volume of the order
of $A(S_1)^{3/2}$ contradicting the hypothesis of the theorem.

\section{Cheeger constants of surfaces}

If $M$ is a riemannian manifold of dimension $n$ one defines the
Cheeger constant $h$ of $M$ by
$$h(M)=\underset {A} { \inf } \, \{\frac {vol_{n-1}(\partial A)}
{vol_n(A)}: vol_n(A)\leq \frac {1} {2} vol_n(M) \}$$ where $A$
ranges over all open subsets of $M$ with smooth boundary. If $M$
is a simplicial manifold one can define the Cheeger constant of
$M$ similarly; now $A$ runs over all simplicial submanifolds of
$M$. To be more precise, we take $A$ to be a union of closed $n$
simplices and we define $\partial A=A\cap \overline {M-A}$. In the
simplicial setting we define $vol_n(A)$ to be the number of
$n$-simplices of $A$ and $vol_{n-1}(\partial A)$ to be the number
of $n-1$-simplices of $\partial A$.

 We remark that this definition
makes sense also if $M$ is more generally an $n$ dimensional
simplicial complex. To make this definition coincide with the
existing literature on graphs one should first take the
barycentric subdivision and then calculate the Cheeger constant.
However as here we are only concerned with surfaces we will not
pass to barycentric subdivisions.

As usual we call $vol_2$ area and $vol_1$ length. If $M$ is a
simplicial 2-manifold or a 2-manifold with a riemannian metric we
denote by $A(M)$ the area of $M$. Similarly if $p$ is a
(simplicial or riemannian) path we denote by $l(p)$ the length of
$p$.

We will show that one can give a bound of the Cheeger constant of
a surface that depends only on its area.
We will treat first the simplicial case and then we will outline
the argument in the riemannian case. In both cases our proof is
based on Besicovitch lemma \cite {Be} (see also \cite{Gr4}, sec.
4.28, p.252, this lemma is sometimes referred to as Almgren's
lemma, \cite {Bo}, \cite{Al}):
\begin{Lem}\label{L:BesicovitchR}  Let $D$ be a riemannian manifold
homeomorphic to the disc and let $\gamma =\partial D$. Suppose
$\gamma $ is split in 4 subpaths, $\gamma =\alpha _1\cup \alpha _2
\cup \alpha _3 \cup \alpha _4$. Let $d_1=d(\alpha _1, \alpha _3)$,
$d_2=d(\alpha _2, \alpha _4)$
 Then
$$A(D)\geq d_1d_2$$
\end{Lem}

We introduce some notation:

If $X$ is a simplicial complex and $K$ is a subcomplex of $X$ we
denote by $star(K)$ the subcomplex of $X$ consisting of all closed
simplices intersecting $K$. We denote by $star_i(K)$ the
subcomplex obtained by repeating the star operation $i$ times. If
$v$ is a vertex of $X$ we define the ball of radius $n$ and center
$v$, $B_v(n)$, by $B_v(n)=star_n(v)$.

 We state below
Besicovitch lemma (\cite {Bo}) in the simplicial setting:
\begin{Lem}\label{L:Besicovitch}  Let $D$ be a simplicial disc and let
$\gamma =\partial D$. Suppose $\gamma $ is split in 4 subpaths,
$\gamma =\alpha _1\cup \alpha _2 \cup \alpha _3 \cup \alpha _4$.
Let $d_1=d(\alpha _1, \alpha _3)$, $d_2=d(\alpha _2, \alpha _4)$
 Then
$$A(D)\geq d_1d_2$$
\end{Lem}
\begin{proof}
We consider $star (\alpha _1)$ and we remark that its boundary has
at least $2d_2$ edges. Since each closed 2-simplex in $star
(\alpha _1)$ intersects the boundary of $star (\alpha _1)$ at at
most 2 edges we conclude that $A(star (\alpha _1))\geq d_2$. Now
we repeat $d_1$ times, i.e. we consider $star _i(\alpha _1)$ for
$i=1,2,...,d_1$ and we remark as before that $$A(star _i(\alpha
_1))\geq d_2$$ It follows that $A(D)\geq d_1d_2$.
\end{proof}
\begin{Rem}
The same inequality applies for disks with a cell complex
structure in which all cells are polygons with 2 or 3 sides
(bigons or triangles). Indeed the proof above applies in this case
too.
\end{Rem}
%
%
%
We start with the inequality for the sphere where the idea of the
proof is more transparent.
\begin{Prop}\label{P:sphere}  Let $S$ be a riemannian manifold
or a simplicial complex homeomorphic to the 2-sphere. Then the
Cheeger constant, $h(S)$, of $S$ satisfies the inequality:
$$ h(S)\leq \frac {16} {\sqrt{A(S)}} $$
where $A(S)$ is the area of $S$.
\end{Prop}
\begin{proof}

We deal first with the simplicial case.
 Let $p$ a closed curve on the 1-skeleton of $S$ of minimal length
dividing $S$ on two regions which have both area bigger or equal
to $\frac {A(S)} {4}$. Let's say $S-p=S_1\sqcup S_2$ (where
$S_1,S_2$ are open). Without loss of generality we assume that
$A(S_1)\geq A(S_2)$. We remark now that there are no `shortcuts'
for $p$ that are contained in $S_1$. More precisely if $a,b\in p$
and $q$ is a path in $S_1$ joining $a,b$ then $l(q)$ is at least
as big as the length of the shortest subpath of $p$ joining $a,b$.
Indeed assume this is not the case. Let's say $p-\{a,b \}=p_1\cup
p_2$ with $l(p_1)\geq l(p_2)$. Without loss of generality we may
assume that $q$ intersects $p$ only at $a,b$. Then $p_1\cup q$ is
a simple closed curve shorter than $p$ which has the same
properties as $p$, a contradiction. We note in particular that
$S_1$ is connected.

We claim that
$$l(p)\leq 4\sqrt {A(S)} $$
This is clearly true if $l(p)<4$.

Otherwise we subdivide $p$ in 4 arcs $p=\alpha _1\cup \alpha _2
\cup \alpha _3 \cup \alpha _4$ such that $l(\alpha _i)\geq \frac
{l(p)}{4}-1$ for all $i$. Since there are no `shortcuts' as we
observed above
$$d(\alpha _1, \alpha _3)\geq \frac {l(p)}{4}-1,\ \ \ d(\alpha _2,
\alpha _4)\geq \frac {l(p)}{4}-1 $$

 Applying lemma \ref{L:Besicovitch} to
$S_1$ we have that $$A(S_1)\geq (l(p)-1)^2/16\Rightarrow 3A(S)\geq
(l(p)-1)^2/4\Rightarrow l(p)\leq 4\sqrt {A(S)}$$

We conclude that $$h(S)\leq \frac {4\sqrt {A(S)}}{A(S)/4}\leq
\frac {16} {\sqrt{A(S)}} $$

We treat now the riemannian case. The argument is along the same
lines. To sidestep the issue of existence and regularity of the
minimal closed curve $p$ we argue instead with $\epsilon $-minimal
curves. More precisely we consider the set $U$ of all simple
closed curves dividing $S$ in two discs which have both area
bigger or equal to $\frac {A(S)} {4}$. Let $L$ be the infimum of
the lengths of the curves in $U$. Given $\epsilon >0$ we consider
$p\in U$ with $l(p)> L-\epsilon$. Let's say that $S-p=S_1\sqcup
S_2$ and $A(S_1)\geq A(S_2)$. Then $p$ does not have $\epsilon
$-shortcuts in $S_1$. That is if $q\subset S_1$ is a path joining
$a,b\in p$ then $l(q)-\epsilon $ is smaller than the length of the
shortest subpath of $p$ joining $a,b$.

We subdivide now $p$ in 4 arcs $p=\alpha _1\cup \alpha _2 \cup
\alpha _3 \cup \alpha _4$ such that $l(\alpha _i)=l(p)/4$ for all
$i$. Since there are no $\epsilon $-`shortcuts' as we observed
above
$$d(\alpha _1, \alpha _3)\geq \frac {l(p)}{4}-\epsilon,\ \ \ d(\alpha _2,
\alpha _4)\geq \frac {l(p)}{4}-\epsilon $$

 Applying lemma \ref{L:BesicovitchR} to
$S_1$ we have that $$A(S_1)\geq (l(p)-\epsilon)^2/16\Rightarrow
12A(S)\geq (l(p)-\epsilon)^2 \Rightarrow l(p)\leq 2\sqrt
{3A(S)}+\epsilon $$

It follows that
$$h(S)\leq \frac {2\sqrt{3A(S)}+\epsilon }{A(S)/4}\Rightarrow h(S)\leq \frac {8\sqrt 3}
{\sqrt{A(S)}} $$ Where the last inequality follows since the
former inequality holds for every $\epsilon >0$. We note that we
obtain a slightly better constant in the riemannian case.

\end{proof}

\begin{Rem}
The same inequality for the Cheeger constant applies for spheres
with a cell complex structure in which all cells are polygons with
2 or 3 sides (bigons or triangles). Indeed the proof above applies
in this case too.
\end{Rem}
%
%
%
%
%

To treat the general case of compact surfaces we need some
technical lemmas.

\begin{Def} Let $S$ be a compact surface with boundary. A simple arc $p$
intersecting the boundary only at its endpoints is said to be
parallel to the boundary if $S-p$ has a contractible component.
Two disjoint simple arcs $p_1,p_2$ intersecting the boundary only
at their endpoints are said to be parallel if $S-(p_1\cup p_2)$
has a contractible component.
\end{Def}

\begin{Lem}\label{L:arcs} Let $S$ be a surface of genus $g$ with
$k$ boundary components. Then there are at most $k+2g$ pairwise
disjoint arcs on $S$ with their endpoints on $\bd S$ such that no
arc is parallel to the boundary and no two arcs are parallel.
\end{Lem}
\begin{proof}
We remark that $\pi _1 (S)$ is a free group of rank $k+2g$. A set
of $i$ arcs which are not parallel pairwise and are not parallel
to the boundary induces a reduced graph of groups decomposition of
$\pi _1 (S)$ with $i$ edges and trivial edge stabilizers. However
the number of edges of any such decomposition can not exceed the
rank of $\pi _1 (S)$.
\end{proof}

\begin{Lem}\label{L:curves} Let $S$ be a closed orientable surface of genus $g\geq 1$.
If $U=\{p_1,...,p_k\}$ is a set of pairwise disjoint simple closed
curves on $S$ such that no component of $S-U$ is contractible then
$k\leq 2g$.
\end{Lem}
\begin{proof}
Without loss of generality we may assume that $U$ is maximal. Then
if we pinch each curve to a point we obtain a space with
fundamental group the free group of rank $2g$, $F_{2g}$. The set
$U$ induces a reduced graph of groups decomposition of $F_{2g}$
with $k$ edges and trivial edge stabilizers, so $k\leq 2g$.

\end{proof}
\begin{Thm}\label{T:cheeger} Let $S$ be a closed orientable
surface of genus $g\geq 1$ equipped either with a riemannian
metric or with a simplicial complex structure. Let $A(S)$ be its
(simplicial or riemannian) area. Then the Cheeger constant,
$h(S)$, of $S$ satisfies the inequality:
$$ h(S)\leq \frac {4\cdot 10^3\cdot g^2} {\sqrt{A(S)}} $$
\end{Thm}
\begin{proof}
We treat first the simplicial case. The proof in the riemannian
case follows the same lines, we outline at the end the changes
which are needed in this case.

 Let $U=\{p_1,...,p_k\}$ be a set of closed
curves on the 1-skeleton of $S$ such that:

1. $S=S_1 \cup S_2$ with $\partial S_1\cap \partial S_2=S_1\cap
S_2=U$, and $A(S_i)\geq A(S)/4,\, i=1,2$.

2. At most one curve $p_i$ bounds a disk in $S$.

3. The sum of the lengths $L=l(p_1)+...+l(p_k)$ is minimal among
all sets of curves satisfying 1,2.

We claim that $$L\leq 10^3\cdot g^2\cdot \sqrt{A(S)}$$

Suppose that this is not the case. By lemma \ref{L:curves}, $k\leq
2g+1$ (note that the $p_i$'s are not necessarily disjoint but can
be made disjoint by pushing them slightly inside $S_1$ or $S_2$).
It follows that there is a curve $p_i\in U$ such that
$$l(p_i)\geq \frac {10^3g^2\sqrt{A(S)}}{2g+1} $$

We set $n=[\sqrt{A(S)} ]+1$. Let's assume that $A(S_1)\geq
A(S_2)$. We remark now that $S_1$ is connected. Indeed suppose
$S_1$ is a disjoint union of two open sets, $T_1,T_2$. Let's say
that $A(T_1)\leq A(T_2)$. We consider $T_1\cup S_2$ and $S_1-T_1$
and we remark that they are separated by a subset of $U$. This
contradicts the minimality of $U$ (property 3). It follows that
$\bd S_1=U$.

We pick now a vertex $v\in p_i$. We claim that $S_1$ is not
contained in the ball of radius $3n$ and center $v$, $B_v(3n)$.
Suppose not. We subdivide $p_i$ at $4g+1$ segments of length
bigger than

$$[\frac {10^3g^2\sqrt{A(S)}}{(2g+1)(4g+1)}]\geq 100n $$
We consider geodesic arcs in $S_1$ from $v$ to the endpoints of
these segments. If some such arc is parallel to the boundary we
can use it to 'cut away' a disc from $S_1$ and contradict the
minimality of $U$ (or in case the disc has area more than half of
the area of $S_1$ we replace $S_1$ by the disc and contradict
property 3).

Otherwise by perturbing these arcs slightly we may arrange so that
they are disjoint. Since we have $4g+1$ arcs by lemma \ref{L:arcs}
two of them are parallel. Using them we can cut away a disk from
$S_1$ (or replace $S_1$ by a disc) which contradicts the
minimality of $U$ (property 3).

We consider now $D_r=B_v(r)\cap S_1$ for $n\leq r\leq 2n$. We
remark that if the length of $\partial D_r\cap S_1$ is bigger than
$2n$ for all $r$ then $A(D_n)\geq n^2\geq A(S)$, a contradiction.

So $\bd D_r$ has length less than $2n$ for some $r$. On the other
hand the length of $D_r\cap p_i$ is bigger than $2n$. So $S_1\cap
D_r$ and $D_r$ have both boundary length smaller than the length
of $\bd S_1$. So we can replace $S_1$ by whichever of the two has
area bigger than $A(S_1)/2$. If this new domain has more than one
boundary component that bounds a disc we just erase this
component. We remark now that the boundary length of the new
domain is smaller than than the boundary length of $S_1$ and this
contradicts the minimality of $U$ (property 3).

The same proof applies in the riemannian case with few changes. We
define a set of closed curves $U=\{ p_1,...,p_k \}$ as before. Now
we may additionally assume that the $p_i$ are simple and disjoint.
To insure this and avoid existence issues we assume that the sum
of their lengths exceeds the minimal possible value by $\epsilon
>0$ (condition 3) among all curves that satisfy 1,2. As before we
argue that there is some $p_i$ such that
$$l(p_i)\geq \frac {10^3g^2\sqrt{A(S)}}{2g+1} -\epsilon $$
We argue as before and w consider $D_r=B_v(t)\cap S_1$ for $r\in
[n,2n]$. Now by the co-area formula if the length of $\partial
D_r\cap S_1$ is bigger than $2n$ for almost all $r$ then
$A(D_n)\geq n^2\geq A(S)$, which gives a contradiction as before.
The rest of the proof applies verbatim to the riemannian case as
well.
\end{proof}

\begin{Rem}\label{R:3dim}  We remark that there is no function of volume that gives an upper bound
for the Cheeger constant of manifolds of dimension higher than 2.
Indeed it's enough to prove this for the ball of dimension 3. We
can obtain examples contradicting the existence of such a bound by
considering sequences of expanders and thickening them.
\end{Rem}

\section{Isoperimetric profiles of surfaces}

If $(M^n,g)$ is a riemannian manifold of infinite volume the
isoperimetric profile function of $M^n$ is a function $I_M:\Bbb R
^+\to \Bbb R ^+$ defined by:
$$ I_M(t)=\underset {\Omega} { \inf } \{vol_{n-1}(\bd \Omega):
\Omega\subset M^n,\, vol_n(\Omega)=t \}$$ where $\Omega $ ranges
over all regions of $M^n$ with smooth boundary. One can define
similarly an isoperimetric profile function $I_M:\Bbb N\to \Bbb N$
for simplicial manifolds $M^n$. In this section we will study
isoperimetric profiles and filling functions of surfaces (so
$vol_2$ is area and $vol_1$ is length).

Other functions related to the isoperimetric problem are the
filling area, $FA_0 $, and homological filling area, $FA^h $,
functions of $M$ that we define now.

If $p$ is a smooth contractible closed curve in $M$ we define its
filling area, $FillA_0(p)$, as follows: We consider all riemannian
discs $D$ such that there is a 1-lipschitz map $f:D\to X$ with
$f|_{\bd D}=p$. We define $FillA_0(p)$ to be the infimum of the
areas of this collection of disks. We define now the filling area
function of $M$ by:

$$ FA_0(t)=\underset {p} { \sup } \{FillA_0(p): l(p)\leq t \}$$

where $p$ ranges over all smooth contractible closed curves of
$M$.

 More generally we can
consider 1-cycles $c$ (i.e. unions of closed curves) that can be
filled by 2-cycles to define the homological filling area
function. To define $FillA(c)$ we consider surfaces with boundary
$(S,\bd S)$, equipped with a riemannian metric, such that there is
a 1-lipschitz map $f:S\to X$ with $f|_{\bd S}=c$. We define then:

$$ FA^h (t)=\underset {c} { \sup } \{FillA(c): l(c)\leq t \}$$
where if $c=c_1\sqcup ...\sqcup c_n$ with $c_i$ closed curves, we
define $l(c)=l(c_1)+...+l(c_n)$. One defines $FA_0$ and $FA^h$
similarly in the simplicial setting as well.

Gromov (\cite {Gr2}, ch. 6, see also \cite{CDP}, ch.6) showed the
following:
\begin{Gromov}\label{T:gromov}  Let $(M^n,g)$ be a simply connected
riemannian manifold. Assume that there is some $t_0$ such that for
all $t>t_0$, $ FA_0(t)\leq \frac {1}{16\pi }t^2$. Then there is a
constant $K$ such that for all $t>t_0$, $ FA_0(t)\leq Kt$.
\end{Gromov}

We remark that Wenger (\cite {We1}) improved $\frac {1}{16\pi }$
to $\frac {1-\epsilon}{4\pi }$ (for any $\epsilon >0$) which is
optimal as the example of the euclidean plane shows. In fact
Gromov's theorem applies more generally to `reasonable' geodesic
metric spaces where a notion of area can be defined (e.g.
simplicial complexes). We note also that Gromov has shown a
stronger (`effective') version than the one we state; it is enough
in fact to have a subquadratic filling for a sufficiently big
range of areas to conclude that the filling is linear.


In the case of surfaces the isoperimetric profile and the filling
area functions are closely related. In fact $FA_0$ is linear for a
space if and only if the space is Gromov hyperbolic (see \cite
{Gr2}). On the other hand if a simply connected surface $S$,
equipped with a riemannian metric, is not Gromov hyperbolic then
for any $t$ there is an embedded loop $\gamma $ in $S$ with
$l(\gamma )>t$ such that $FillA_0(\gamma )> \frac {1}{16 \pi}
l(\gamma )^2$. If $t_1=FillA_0(\gamma )$ we see that for any $t>0$
there is some $t_1>t$ such that $I_S(t_1) < 4 \sqrt \pi \sqrt
t_1$.

\begin{Thm}\label{T:gromov1}  Let $(S,g)$ be a simply connected
riemannian surface. Assume that there is some $t_0$ such that for
all $t>t_0$, $ I_S(t)\geq 4\sqrt \pi \sqrt t$. Then there is a
constant $\delta >0 $ such that for all $t>t_0$, $ I_S(t)\geq
\delta t$.
\end{Thm}

We remark that in many cases $FA _0$ and $FA ^h$ are equal (e.g.
this holds for the Euclidean and Hyperbolic plane).

This does not hold always however. We give now some examples to
illustrate the relationship between the filling area functions and
the isoperimetric profile. If $f(t),g(t)$ are functions we write
$f(t)\sim g(t)$ if $$\limsup _{t\to \infty } \frac {f(t)}{g(t)}<
\infty, \ \ \  \liminf _{t\to \infty } \frac {f(t)}{g(t)}>0 $$

\begin{Ex} Let $X$ be the punctured Euclidean plane $\Bbb R^2- \Bbb
Z^2$. Then $FA_0(t)\sim t$ (see \cite{BE}, \cite{PS}) while
$FA^h(t)\sim t^2$. The isoperimetric profile $I_X$ is the inverse
of $FA^h$ so $I_X(t) \sim \sqrt t$.
\end{Ex}

\begin{Ex} Let $X$ be the cylinder $S^1\times \Bbb R$ with the
standard product metric. Then $FA_0(t)\sim t^2 $. Indeed for any
$X$, $FA_0$ is the same for $X$ and for the universal covering
$\tilde X$. Here $\tilde X=\Bbb E ^2$ (the euclidean plane).  On
the other hand if $s$ is the length of the $\Bbb S^1$ factor for
any $t>2s$ we have $FA^h(t)=\infty $. Similarly for the
isoperimetric profile there is some $t_0$ such that for all
$t>t_0$, $I_X(t)=2s$.
\end{Ex}

\begin{Ex} Let $X_n$ be the space obtained by removed a ball of
radius $n$ from the hyperbolic plane $\Bbb H ^2$. Let $X$ be
isometric to the hyperbolic plane. We fix a point $O\in X$ and we
consider a sequence of points $x_n$ such that $d(x_n,O)=2^n$. For
each $n$ we remove the disk of radius $n$ and center $x_n$ from
$X$ and we glue along the boundary of the disk a copy of $X_n$.
The space $Y$ obtained has $FA_0(t)\sim FA^h(t) \sim I_Y(t) \sim
t$. We remark that $Y$ is not Gromov hyperbolic.
\end{Ex}
%

%
%

It is reasonable to ask whether Gromov's theorem extends to all
 surfaces. The answer is no in general
but we can show that the theorem holds for surfaces of bounded
genus (this applies for example to riemannian planes with
infinitely many holes, compare \cite{PS}). We have the following:

\begin{Thm}\label{T:profile1}  Let $S$ be a
plane with holes equipped either with riemannian metric or with a
simplicial complex structure. Assume that there is some $K>0$ such
that for all $t\in [K,100K]$, $ I_S(t)\geq 10^2 \sqrt t$. Then
there is a constant $\delta >0 $ such that for all $t>K$, $
I_S(t)\geq \delta t$.
\end{Thm}
\begin{proof}
We treat the simplicial case first.
 We will show that $I(t)\geq \frac {1}{\sqrt K}t$ for all $t>K$.

%
%


We argue by contradiction. So let $c$ be a 1-cycle with minimal
filling area and $FillA(c)>K$ such that $l(c)< \frac {1}{\sqrt
K}FillA(c)$. Let's say that $c=\bd R$. We claim that
$FillA(c)>100K$. Indeed if $FillA(c)\leq 100K$ then $l(c)\geq
100\sqrt {FillA(c)}$ hence
$$\frac {1}{\sqrt K}FillA(c)> 100\sqrt {FillA(c)}\Rightarrow
FillA(c)>10^4K$$

By our minimality assumption $R$ is connected, so $R$ is a sphere
with holes (possibly a disc). We collapse all holes to points and
we obtain a sphere $\Sigma $ will a cell complex structure in
which all cells are either bigons or triangles.


In this way we obtain a map $f: R\to \Sigma $  which is 1-1 on
open 2-simplices from $R$.

By Proposition \ref{P:sphere} (and the remark following it) there
is a simple closed curve $p$ in $\Sigma ^1$ such that
$FillA(p)\geq A(\Sigma )/4$ and
$$l(p)\leq 4\sqrt {A(\Sigma )} $$

It follows that $$FillA(p)\geq \frac {1}{64} l(p)^2$$

We consider now the set of curves $q$ in $\Sigma ^1$ with filling
area $FillA(q)\geq K$ that satisfy $FillA(q)> \frac {1}{100}
l(q)^2$. Clearly this set is not empty. Let $\gamma $ in $\Sigma
^1$ of minimal filling area with this property. We will show that
$4K\geq FillA(\gamma )\geq K$. Assume this is not the case.

We subdivide $\gamma $ in 4 arcs $\gamma =\alpha _1\cup \alpha _2
\cup \alpha _3 \cup \alpha _4$ such that $l(\alpha _i)\geq \frac
{l(\gamma )}{4}-1$ for all $i$.

We claim that
$$d(\alpha _1, \alpha _3)\geq \frac {l(\gamma )}{4},\ \ \ d(\alpha _2,
\alpha _4)\geq \frac {l(\gamma )}{4} $$

We argue by contradiction. Assume that $d(\alpha _1, \alpha _3)<
\frac {l(\gamma )}{4}$ and let $w$ be a path from $\alpha _1$  to
$\alpha _3$ of length $r\leq \frac {l(\gamma )}{4}$. Then, using
$w$, we split $\gamma $ into two curves $\gamma _1, \gamma _2$
such that $\gamma _1\cap \gamma _2=w$, $\gamma _1\cup \gamma
_2=\gamma \cup w$. Let's say $l(\gamma_1)\geq l(\gamma _2)$. Then
$$l(\gamma _1)=\frac {l(\gamma )}{2}+a+r,\ \  l(\gamma _2)=\frac
{l(\gamma )}{2}-a+r $$

for some $a<l(\gamma )/4$. To simplify the notation we set
$l=l(\gamma )$. Since we assume that $FillA(\gamma )>2K$ we have
that $FillA(\gamma _1)>K, \, FillA(\gamma _2)>K$, so

$$FillA(\gamma )\leq FillA(\gamma _1)+ FillA(\gamma _2)\leq \frac {1}{100}[ (\frac
{l}{2}+a+r)^2+(\frac {l}{2}-a+r)^2]=$$ $$=\frac {1}{100}(\frac
{l^2}{2}+2a^2+r^2+lr)$$

Since $a\leq \frac {l}{4}, r\leq \frac {l}{4}$ we have

$$\frac{l^2}{2}+2a^2+r^2+lr\leq
\frac{l^2}{2}+\frac{l^2}{8}+\frac{l^2}{16}+\frac{l^2}{4}<l^2$$

which is a contradiction. We may now apply Lemma
\ref{L:Besicovitch} to $\gamma $ and conclude that $FillA(\gamma
)\geq \frac {l(\gamma )^2}{16} $ which is again a contradiction.
We conclude that $4K\geq FillA(\gamma )\geq K$.

We lift now $\gamma $ via $f$ to $R$. $\gamma $ lifts to a set of
arcs (or a single simple closed curve) that separate $R$ into two
2-chains $R_1,R_2$. Let's denote this set of arcs by $\alpha $.
Let's say that $K\leq A(R_1)\leq 4K$. Then $\partial R_1=c_1\cup
\alpha $ and $\partial R_2=c_2\cup \alpha $ with $c_1\cup c_2=c$.
By our assumption on $c$ we have
$$l(c_1)+l(c_2)\leq \frac{1}{\sqrt K}(A(R_1)+A(R_2))\ \ \  (1)$$
On the other hand since $c$ is minimal with this property we have
$$l(c_1)+l(\gamma )\geq \frac{1}{\sqrt K}A(R_1)\ \ \ (2)$$
$$l(c_2)+l(\gamma )\geq \frac{1}{\sqrt K}A(R_2)\ \ \ (3)$$
By the way $\gamma $ was defined we have
$$A(R_1)\geq \frac {1}{100} l(\gamma )^2$$
Since $A(R_1)\leq 4K$ we have $l(\gamma )\leq 20 \sqrt K$. From
the hypothesis of the theorem since $A(R_1)\in [K,100K]$ we have

$$ 100A(R_1)\leq  l(\gamma )+l(c_1)\Rightarrow
l(c_1)\geq 80 \sqrt K $$

Substituting in (1) we obtain
$$l(c_2)\leq \frac{1}{\sqrt K}(4K+A(R_2))-80\sqrt K$$
Therefore
$$l(c_2)+l(\gamma )\leq \frac{1}{\sqrt K}A(R_2)-76\sqrt K$$
and from (3)
$$\frac{1}{\sqrt K}A(R_2)\leq \frac{1}{\sqrt K}A(R_2)-76\sqrt K$$
which is a contradiction.

The proof in the riemannian case is identical. One has just to
note that when we collapse the boundary curves to points we obtain
a riemann metric with some singularities. Our estimates for
Cheeger constants apply however to this case as well. One can see
this e.g. by approximating the singular metric by a non singular
one or by noting that our proof of the Cheeger constant bounds
work also for singular metrics.
\end{proof}

\begin{Cor}\label{C:profile2}  Let $S$ be a non-compact
surface of finite genus equipped either with a riemannian metric
or with a simplicial complex structure. Assume that there is some
$K>0$ such that for all $t\in [K,100K]$, $ I_S(t)\geq 10^2 \sqrt
t$. Then there is a constant $\delta >0 $ such that for all $t>K$,
$ I_S(t)\geq \delta t$.
\end{Cor}
\begin{proof}
There is a finite set $\{p_1,...,p_n\}$ of smooth, rectifiable,
simple closed curves (or a finite set of simple closed curves
lying in $S^1$ is the simplicial case) such that
$$S-\{p_1,...,p_n\}=B\sqcup B_1\sqcup....\sqcup B_k$$
with $B$ a surface of finite area and $B_1,...,B_k$ planes with
holes. From the previous theorem we have that there are $\delta
_1,...,\delta _k$ such that for all $t>K$ we have
$$I_{B_i}(t)\geq \delta _i t,\ \, i=1,...,k$$

Let $\delta '= \min \{\delta _1,...,\delta _k \}$. We set $V=A(B)$
and $L=l(p_1)+...+l(p_n)$.

If $\Omega $ is a domain in $S$ with rectifiable boundary there is
some $i\in \{1,...,k\}$ such that $$A(\Omega \cap B_i)\geq \frac
{A(\Omega )-V}{k}$$

Let $c=\bd \Omega \cap B_i$. Since $\bd (\Omega\cap B_i)\subset
c\cup \{p_1,...,p_n\}$ we have that $l(\bd (\Omega\cap B_i))\leq
L+l(c)$. If

 $$\frac{A(\Omega )-V}{k}>K\Leftrightarrow A(\Omega )>kK+V$$ we have
$$L+l(c)\geq \delta '\frac {A(\Omega )-V}{k}\Rightarrow l(c)\geq \delta '\frac {A(\Omega )-V}{k}-L$$
It follows that if $$\delta '\frac {A(\Omega )}{k}\geq 2(\frac {
\delta ' V}{k}+L)\Leftrightarrow A(\Omega)\geq 2V+\frac
{2kL}{\delta '}$$ we have
$$l(\bd \Omega )\geq l(c)\geq \delta '\frac {A(\Omega )}{2k} $$
We conclude that for all $$t>\max (kK+V,2V+\frac {2kL}{\delta
'})$$ we have
$$I_S(t)\geq \frac {\delta '}{2k}t$$

We note further that if $\delta _1=\inf \{I_S(t):t>K \}$ and
$$\delta =\min (\frac {\delta _1}{K}, \frac {\delta '}{2k})$$ then
for all $t>K$ we have
$$I_S(t)\geq  \delta t$$

\end{proof}

\begin{Rem}\label{R:FA} The previous theorem implies that if the
filling area function is subquadratic for surfaces of finite genus
then it is actually linear. In fact one may give a similar proof
to another generalization of Gromov's theorem. Let $X$ be either a
riemannian manifold or a simplicial complex. Let $c$ be a 1-cycle.
If $c=\bd R$ for some 2-chain $R$ then we define the genus of $R$
to be the genus of the 2-cycle we obtain from $R$ by collapsing
$c$ to a point. If $c$ is a 1-cycle in $X$ we define the
$g$-filling area of $c$ by
$$FillA_g(c)=\inf \{A(S): S \text{ is a 2-chain of genus at most
$g$ s.t. } \bd S=c \}$$ Note that with this definition $FillA_0$
is slightly more general than before as it applies to 1-cycles and
not just closed curves but this does not affect what follows.
 We
define now the $g$-filling area of $X$ by
$$FA_g(t)=\underset {c } \sup \{ FillA_g(c):  c=\bd R, \text {  $R$ of genus $\leq g $ and } l(c)\leq t \}$$

With this notation Gromov's theorem says that if $FA_0$ is
subquadratic then it is bounded by a linear function. In fact now
we can generalize this for any $g$: For any $g$ there is some
$\epsilon _g>0$ such that if for some $t_0$, $FA_g(t)\leq
\epsilon_g t$ for all $t>t_0$ then there is some $K>0$ such that
$FA_g(t)\leq K t$ for all $t>t_0$. The proof goes along the same
lines as the proof of theorem \ref{T:profile1}. We argue by
contradiction assuming that we have a minimal $1$-cycle $c$ that
violates the linear isoperimetric inequality. We fill it by a
minimal area 2-chain $S$ of genus at most $g$. Then we collapse
all boundary components of $S$ to obtain a closed 2-cycle $\Sigma
$ of genus $\leq g$. Using theorem \ref{T:cheeger} we show that we
can `cut' a 2-chain from $\Sigma $ with small boundary length and
big area. We lift this back to $S$ and we argue as in theorem
\ref{T:profile1} to contradict our assumption that $c$ is minimal.

\end{Rem}

\section{Isoperimetric inequalities}

In this section we will study the question whether Gromov's `gap'
theorem for $FA_0$ extends to the 3-dimensional filling function
$FV_3$. The filling area function $FA_0$ is important for group
theory since it is related to the word problem. In fact if $G$ is
a finitely presented group and $X$ is its Cayley complex then $G$
has a solvable word problem if and only if $FA_0(X)$ is bounded by
a recursive function ($X$ might not be a simplicial complex but
one can pass to a simplicial subdivision to make sense of
$FA_0(X)$). The question whether there are other `gaps' for $FA_0$
apart from between $n$ and $n^2$ for finitely presented groups was
answered in the negative (see \cite{SBR}, \cite{Br}, \cite{BB}).
It is easy to see that one can produce simplicial (or riemannian)
planes with $FA_0$ of the form, say, $n^r,\,r\in (2,\infty )$ and
Grimaldi-Pansu (\cite{GP}) study the finer question of
characterizing completely filling functions for riemannian planes.

Gromov (\cite{Gr3}) has given estimates and formulated conjectures
for higher dimensional filling functions of nilpotent groups (see
also \cite{WP} and \cite{BBFS} for interesting examples of higher
filling functions of groups).

We note that as we move to higher dimensions we have two possible
ways to define filling functions. We can either define them by
considering fillings of (singular) spheres by balls or more
generally one may consider filling of higher dimensional cycles
(e.g. filling of orientable surfaces of genus $g\geq 0$ in
dimension 3). Here we take the second option, apart from being
easier to define technically it seems more natural. For example,
as it is shown in \cite{P2} filling of 2-spheres in groups is
always subrecursive in contrast to $FA_0$ which is not
subrecursive for groups with unsolvable word problem. So examining
filling only of 2-spheres seems quite restrictive.

We refer to the introduction for the definition of the terms in
the theorem below. To simplify notation we denote $Fillvol_2(c)$
by $FillA(c)$ if $c$ is a 1-cycle and $Fillvol_3(S)$ by $FillV(S)$
if $S$ is a 2-cycle.

\begin{Thm}\label{T:dim3} Let $X$ be a simplicial
complex such that $H_1(X)=H_2(X)=0$. Assume that the following
hold:
\begin{itemize} \item There is some $K>0$ such that
 $FV_2(n)\leq Kn^2$ for all $n\in \Bbb N$.

 \item $$\lim _{n\to \infty }\frac {FV_3(n)}{n^{\frac {3}{2}}}=0$$

 \item There is some $g\in \Bbb N$ such that if $S$ is a minimizer
 $2-cycle $ in $X^2$ then $S$ is represented by a surface of genus at
 most $g$.

\end{itemize}
Then for every $\epsilon >0$ we have
$$\lim _{n\to \infty }\frac {FV_3(n)}{n^{1+\epsilon }}=0$$

\end{Thm}

We don't know whether the third condition on the bound of the
genus of the minimizers is in fact necessary. It is quite crucial
however for our proof. We use it to deduce that the diameter of a
minimizer $S$ is of the order of $\sqrt {A(S)}$. This in turn is
based on our upper bound of Cheeger constants in terms of genus.
So our proof would work as well if we assumed that there is an
upper bound for the Cheeger constant of minimizers of the form
$\frac {c } {\sqrt {A(S)}}$ or if we assumed that there is a lower
bound for the diameter of minimizers of the form $\delta  \sqrt
{A(S)}$. It would be interesting to remove the condition on
minimizers even in the case that $X$ is a non-positively manifold
homeomorphic to $\Bbb R ^3$. We remark that it is not known
whether the isoperimetric profile of a non-positively curved
manifold homeomorphic to $\Bbb R ^n$ is dominated by the
isoperimetric profile of the Euclidean space $\Bbb E ^n$ (this is
known however for $n\leq 4$, see \cite {Cr}, \cite{Kl}).

We are going to prove a somewhat stronger statement that implies
theorem \ref{T:dim3}:
\begin{Thm}\label{T:dim3eff} Let $X$ be a simplicial
complex such that $H_1(X)=H_2(X)=0$ Assume that the following
hold:
\begin{itemize} \item There is some $K>0$ such that
 $FV_2(n)\leq Kn^2$ for all $n\in \Bbb N$.

 \item There is some $g\in \Bbb N$ such that if $S$ is a minimizer
 $2-cycle $ in $X^2$ then $S$ is represented by a surface of genus at
 most $g$.

\end{itemize}

Then given $\epsilon >0$ there is a constant $\alpha =\alpha
(K,g)>0$ such that the following holds: If there is an $n_0$ such
that for all $n>n_0$, $FV_3(n)<\alpha n^{\frac {3}{2}}$, then
$$\lim _{n\to \infty }\frac {FV_3(n)}{n^{1+\epsilon }}=0$$
\end{Thm}
\begin{proof}
 In the course of the proof we will
need to introduce some new constants; we will indicate the
previous constants that the new constants depend on, e.g. for the
new constant $c$ we write $c=c(A,B)$ to indicate that $c$ depends
on the previous defined constants, $A,B$. It is possible always to
give explicit estimates for the constants but we refrain from
doing this as we don't find it instructive.

To show the theorem it is enough to show that for any $\epsilon
>0$ there is some $\beta =\beta (K,g)>0$ with the following property: If there is an $n_0$ such
that for all $n>n_0$, $FV_3(n)<\beta n^{\frac {3}{2}}$ then
$$\limsup_{n\to \infty }\frac {FV_3(n)}{n^{1+\epsilon }}\ne \infty $$
Indeed we can then take $\alpha (\epsilon )=\beta (\epsilon/2)$.

We argue by contradiction. The value of $\beta $ will be specified
in the course of the proof.
 So we assume that for some $\epsilon
>0$ the following holds:

For any $M>0$ there is some 2-cycle $S$ such that $FillV(S)>
MA(S)^{1+\epsilon }$.

\begin{Lem}\label{L:diameter} There is a $\delta =\delta (\epsilon,g, K)>0$ such that for any $M>0$ if
$S$ is a 2-cycle of minimal area such that $FillV(S)>
MA(S)^{1+\epsilon }$ then $diam (S)>\delta \sqrt {A(S)}$.

\end{Lem}
\begin{proof}

We set $n=\sqrt {A(S)}$. As we saw in the proof of theorem
\ref{T:cheeger} there is a decomposition of $S$ in two pieces
$S_1,S_2$ such that: \begin{itemize}

\item $S=S_1\cup S_2$
\item $S_1\cup S_2=\partial S_1\cap \partial S_2$
\item $A(S_i)\geq A(S)/4,\, i=1,2 $
\item $l(S_1\cap S_2)\leq 10^3g^2n$
\item $S_1\cap S_2$ has at most $2g+1$ components.
\end{itemize}

We claim now that if $p$ is a closed curve on $X^1$ of diameter
less than $\delta l(p)$ then $FillA(p)\leq 40K\delta l(p)^2$.

To see this subdivide $p$ into $[1/\delta ]+1$ segments of length
at most $\delta l(p)+1$. Let $v_1,...,v_r$ be the successive
endpoints of these segments. We consider geodesic segments
$[v_1,v_3],...,[v_1,v_{r-1}]$ and we use them to break $p$ into
$r-2$ loops each of which has length at most $3\delta l(p)+3$.
Since $FillA(p)$ is less or equal to the sum of the areas of these
loops we have:
$$FillA(p)\leq K([1/ \delta ]-1)(3\delta l(p)+3)^2\leq K\frac {1}{\delta }(6\delta l(p))^2\leq 40K\delta l(p)^2$$

Assume now that $diam (S)\leq \delta n$. We will show that this
leads to a contradiction if $\delta $ is too small.

Let $c=S_1\cap S_2$. Since $c$ has at most $2g+1$ components and
$diam(c)\leq diam (S)\leq \delta n$ using the above estimate for a
single simple closed curve we obtain for $c$:
$$FillA(c)\leq 40K(2g+1)\delta l(c)^2\leq 40K(2g+1)\delta (10^3g^2n)^2$$
Let $\tilde S$ be a 2-cycle filling $c$ with $A(\tilde S)\leq
40K(2g+1)\delta l(c)^2$. We break $S$ into two 2-cycles using
$\tilde S$: $\tilde S_1=S_1+\tilde S$ and $\tilde S_2=S_2-\tilde
S$. We set $\delta '=40K(2g+1)10^6g^4\delta $. So $$A(\tilde
S)\leq \delta ' A(S)$$ If $\delta $ is sufficiently small
$A(\tilde S)$ is smaller than $A(S_1), \, A(S_2)$. Using the
minimality of $S$ we have:

$$FillV(S)\leq FillV(\tilde S_1)+FillV(\tilde S_2)\leq MA(\tilde S_1)^{1+\epsilon }+MA(\tilde S_2)^{1+\epsilon
}$$

Since $A(S_1)$ and $A(S_2)$ are bigger than $A(S)/4$ there is some
$a\in [\frac {1}{4},\frac {3}{4}]$ such that $A(S_1)=aA(S)$ and
$A(S_2)=(1-a)A(S)$. We have $$A(\tilde S_1)\leq aA(S)+\delta '
A(S),\,\, A(\tilde S_2)\leq aA(S)+\delta 'A(S)$$

Substituting above we have

$$FillV(S)\leq MA(S)^{1+\epsilon }[(a+\delta ')^{1+\epsilon
}+(1-a+\delta ')^{1+\epsilon }] $$

Since the function $x^{1+\epsilon }$ is strictly convex
$a^{1+\epsilon }+(1-a)^{1+\epsilon}<1$ for all $a\in [1/4,3/4]$.
It follows that if $\delta '$ is small enough $(a+\delta
')^{1+\epsilon }+(1-a+\delta ')^{1+\epsilon }<1$. Clearly one can
give an explicit estimate for $\delta '$ in terms of $\epsilon $.

Now if $$\delta \leq \frac {\delta '}{40K(2g+1)10^6g^4}$$ we have

  $$FillV(S)\leq MA(S)^{1+\epsilon
}$$ which is a contradiction.

\end{proof}

We need a technical lemma:
\begin{Lem}\label{L:ineq}

Given $\epsilon >0$ there is some $\lambda >0$ such that for any
$x\in (0,1/2]$ the following inequality holds:
$$(x+\lambda x)^{1+\epsilon }+(1-x+\lambda x)^{1+\epsilon }<1 $$

\end{Lem}
\begin{proof}

We consider the function $$f(x)=1-(x+\lambda x)^{1+\epsilon
}-(1-x+\lambda x)^{1+\epsilon }$$

We have $$f'(x)=(1+\epsilon)[-(1+\lambda )^{1+\epsilon
}x^{\epsilon }+(1-\lambda )(1-x+\lambda x)^{\epsilon }]$$

We remark now that there is a constant $c>0$ such that if $\lambda
<1/2$ we have
$$-(1+\lambda )^{1+\epsilon
}x^{\epsilon }+(1-\lambda )(1-x+\lambda x)^{\epsilon }>0,\,
\forall x\in [0,c]$$ Since $f(0)=0$ we conclude that $f(x)>0$ for
$x\in [0,c]$, if $\lambda <1/2$.

Now we remark that the function $x^{1+\epsilon }$ is strictly
convex. It follows that $x^{1+\epsilon }+(1-x)^{1+\epsilon }$
restricted on the interval $[c,1/2]$ is strictly smaller than 1.
It follows that there is some $\lambda >0$ such that $$(x+\lambda
x)^{1+\epsilon }+(1-x+\lambda x)^{1+\epsilon }<1$$ for all $x\in
[c,1/2]$. So there is some $1/2>\lambda >0$ such that for any
$x\in (0,1/2)$ we have
$$(x+\lambda x)^{1+\epsilon }+(1-x+\lambda x)^{1+\epsilon }<1$$
\end{proof}

 In what follows given $M>0$ we consider
a 2-cycle $S$ of minimal area such that $FillV(S)>
MA(S)^{1+\epsilon }$. Let $R$ be a 3-chain such that $\bd R=S$ and
$V(R)=FillV(S)$. We consider $R$ as a subset of $X$. We fix a
vertex $x\in S$ and we consider $B_i(x)$ in $X$. Let $R=\sum
n_j\sigma _j$ with $n_j\in \Bbb Z$. We define $R_i$ to be the
chain:

$$R_i=\underset {\sigma _k\in B_i(x)} \sum n_k\sigma _k $$

We consider now all decompositions of $\bd R_i$ as sum of two
chains $\bd R_i=R_1+R_2$. We consider the minimal value of
$A(S-R_1)$ over all such decompositions. Let $\bd R_i=O_i+I_i$ be
a decomposition of $\bd R_i$ such that $A(S-O_i)$ attains this
minimum.

 With this
notation we have the following lemma.
\begin{Lem}\label{L:outer} There is a $\lambda >0$ such that for
any $M>0$ if $S$ is a 2-cycle $S$ of minimal area such that
$FillV(S)\geq MA(S)^{1+\epsilon }$ then the following holds:

$$A(I_i)\geq \lambda \min \{A(O_i),A(S-O_i) \}, \ \  \ \ \forall i\in \textstyle [\frac{\delta
n}{4},\frac{\delta n}{2}]$$ where $\delta $ is given by lemma
\ref{L:diameter}.

\end{Lem}
\begin{proof}
Let $\lambda $ be as in lemma \ref{L:ineq}. We argue by
contradiction, ie we assume that the inequality of the lemma does
not hold for some $i$. We consider the 2-cycles $\bd R_i$ and
$S-\bd R_i$. We remark that $A(\bd R_i)$ and $A(S-\bd R_i)$ are
both smaller than $A(S)$. By our assumption on $S$ we have the
inequalities:

$$FillV(\bd R_i)< MA(\bd R_i)^{1+\epsilon }$$

$$FillV(S-\bd R_i)< MA(S-\bd R_i)^{1+\epsilon }$$

We also have

$$ FillV(S)\leq FillV(\bd R_i)+FillV(S-\bd R_i)$$

Hence

$$ FillV(S)< M[A(\bd R_i)^{1+\epsilon }+A(S-\bd R_i)^{1+\epsilon }] \ \ \ \ (1)$$

Now $A(\bd R_i)=A(O_i)+A(I_i)$ and $A(S-\bd R_i)=A(I_i)+A(S-O_i)$.
Let $$m=max\{A(O_i),A(S-O_i)\}$$ We set $$x=\frac {m}{A(S)}$$ From
inequality (1) and from our assumption we obtain:

$$  FillV(S)< MA(S)^{1+\epsilon }[(x+\lambda x)^{1+\epsilon}+(1-x+\lambda
x)^{1+\epsilon}]< MA(S)^{1+\epsilon } $$

where the last inequality follows from lemma \ref{L:ineq}. This is
clearly a contradiction.

\end{proof}

\begin{Lem}\label{L:quadratic}  Let $\mu >0$, $a>1$ and let $f:[a,2a]\to \Bbb R ^+$ be
a continuous function such that $a<f(t)<\mu a^2$ for all $t\in
[a,2a]$. If
$$F(s)=\int _a ^s f(t)dt $$
then for some $r\in [a,2a]$ we have $$F(r)> \frac {1}{\sqrt{3\mu
}}f(r)^{\frac{3}{2}}$$
\end{Lem}
\begin{proof}
We consider the function $$g(t)=\frac {\mu (t-a)^2}{3}$$ Clearly
$f(0)>g(0)$ and $g(2a)>f(2a)$. Let $r$ be the first point in
$[a,2a]$ such that $g(r)=f(r)$. Then
$$F(r)=\int _a ^r f(t)dt > \int _a ^r \frac {\mu (t-a)^2}{3}dt=\frac
{\mu (r-a)^3}{9}=\frac {1}{\sqrt{3\mu }}f(r)^{\frac{3}{2}}$$

\end{proof}

We will need a `discrete' version of the above lemma, which we
state now:
\begin{Lem}\label{L:dquad}  Let $\mu >0$, $a\in \Bbb N$ and let $f(i)\in \Bbb R ^+$ be
a sequence such that $a<f(i)<\mu a^2$ for all $i\in \Bbb N\cap
[a,2a]$. If
$$F(s)=\sum _{i=a} ^s f(i) $$
then for some $r$ we have $$F(r)> \frac {1}{\sqrt{3\mu
}}f(r)^{\frac{3}{2}}$$
\end{Lem}
\begin{proof}
We extend $f$ to a piecewise constant function defined on
$(a-1,2a]$ by posing $f(t)=f(i)$ for all $t\in (i-1,i)$ where
$i=a,a+1,...,2a$.

 We consider the function $$g(t)=\frac {\mu (t-a)^2}{3}$$
Clearly $f(0)>g(0)$ and $g(2a)>f(2a)$.

 Let $x=\inf \{t\in [a,2a]:g(x)\geq f(x)\} $. Then
$$\int _a ^x f(t)dt > \int _a ^x \frac {\mu (t-a)^2}{3}dt=\frac
{\mu (x-a)^3}{9}$$    If $x\in \Bbb N$ and $g(x)=f(x)$ we take
$r=x$ and we have
$$\sum _{i=a} ^r f(i)>\frac
{\mu (x-a)^3}{9}=\frac {1}{\sqrt{3\mu }}f(r)^{\frac{3}{2}}$$

If $x\in \Bbb N$ and $g(x)<f(x)$ we remark that $f(x+1)<f(x)$ so
the desired inequality holds for $r=x+1$. Otherwise let $r=[x]+1$.
We have $g(x)=f(x)=f(r)$ and $2g(x)\geq g(r)$ so we have:
$$F(r)>\int _a ^x f(t)dt >\int _a ^x \frac {\mu (t-a)^2}{3}dt=\frac {1}{\sqrt{3\mu
}}f(x)^{\frac{3}{2}}=\frac {1}{\sqrt{3\mu }}f(r)^{\frac{3}{2}}  $$

\end{proof}

We take now $M=M(n_0)$ `sufficiently big' and we consider a
2-cycle $S$ of minimal area such that $FillV(S)> MA(S)^{1+\epsilon
}$. We will explain how we choose $M$ at the relevant point of the
proof.

 Let $R$ be a 3-chain such that $\bd R=S$ and $V(R)=FillV(S)$.

Let $\delta $ be as in lemma \ref{L:diameter}. We fix a vertex
$x\in S$ and we consider $B_i(x)\cap S$ and $B_i(x)\cap R$ for
$\delta n/4\leq i\leq \delta n /2$. We define $R_i$, $O_i$, $I_i$
as above. We remark now that the following inequalities hold:

$$A(O_i)\geq 2i\geq \delta n /2$$

$$A(S-O_i)\geq \delta n /2$$

$$V(R_i-R_{i-1})\geq A(I_{i})/3\geq \lambda \delta n /2$$

where $\lambda $ is the constant provided by lemma \ref{L:outer}
and the last inequality follows from the same lemma.

We consider now the finite sequence $$f(i)=A(I_i),\, i\in [\delta
n/4, \delta n/2]$$ We set $a=\delta n/4$ and we remark that
$f(i)\geq a$ for all $i\in [\delta n/4, \delta n/2]$ and $f(i)\leq
4n^2\leq \frac {64}{\delta ^2}a^2$. Applying lemma \ref{L:dquad}
to $f$ with $\mu = \frac {64}{\delta ^2}$ we conclude that there
is some $r\in [\delta n/4, \delta n/2]$ such that

$$\sum _{i=a} ^r f(i) > \frac {1}{\sqrt{3\mu
}}f(r)^{\frac{3}{2}}$$

We have also the inequality:

$$V(R_r)\geq \sum_{j=\delta n/4 }^r f(j)/3\geq \frac {1}{3\sqrt{3\mu
}}f(r)^{\frac{3}{2}}$$

From lemma \ref{L:outer} $$f(r)\geq \lambda \min \{A(O_r),A(S-O_r)
\}$$

We distinguish now two cases:

Case 1. $f(r)\geq \lambda A(O_r)$. In this case
$$A(\bd R_r)=A(O_r)+A(I_r)\leq (1+\frac {1}{\lambda })f(r)$$
while $$V(R_r)\geq \frac {1}{3\sqrt{3\mu }}f(r)^{\frac{3}{2}} $$

So if $$\beta = \frac {1}{3\sqrt{3\mu }(1+\frac {1}{\lambda
})^{\frac{3}{2}} }$$

and $M$ is big enough so that $A(\bd R_r )>n_0$ we have

$$V(R_r)\geq \beta A(\bd R_r)^{\frac{3}{2}} $$
which contradicts the our assumption.

Case 2. $f(r)< \lambda A(O_r)$. Then $f(r)\geq A(S-O_r)$ so
$A(S-O_r)\leq A(O_r)$. In this case we pick $y\in S$ with
$d(x,y)=diam S^1$ and we repeat the construction considering
$B_y(i)$ instead of $B_x(i)$. We obtain a 3-chain as before which
we denote $R'_{r'}$. $\bd R'_{r'}\cap S\subset A(S-O_r)$ we obtain
a contradiction from $\bd R'_{r'}$ as in case 1.

\end{proof}

\begin{Rem}
The assumption that $FV_2(n)$ is bounded by a quadratic function
does not play an essential role in the proof above. One may
substitute this by $FV_2(n)\leq Kn^r$ for some $K>0,r>2$ and
change the conclusion to: If $$\lim _{n\to \infty } \frac
{FV_3(n)}{n^{\frac {2r-1}{2r-2}}}=0 $$ then $$\lim _{n\to \infty }
\frac {FV_3(n)}{n^{1+\epsilon }}=0 $$ for any $\epsilon >0$ . This
shows that there is some relationship between $FV_2(n)$ and
$FV_3(n)$, always of course under the assumption of the bound on
the genus of minimizers.

\end{Rem}

\end{document}
\bye